\documentclass[11pt]{article}
\usepackage[margin=1in]{geometry}

\usepackage{amsmath} %
\usepackage{amssymb}  %
\usepackage{amsthm}
\usepackage{hyperref}
\usepackage{algorithm}
\usepackage{algpseudocode}

\theoremstyle{plain}  %
\newtheorem{theorem}{Theorem}[section]
\newtheorem{lemma}[theorem]{Lemma}

\newtheorem{definition}[theorem]{Definition}

\theoremstyle{definition}

\theoremstyle{remark} %
\newtheorem{example}{Example}[section]
\newtheorem{remark}{Remark}[section]

\title{\LARGE \bf A convex polynomial that is not sos-convex}
\author{Amir Ali Ahmadi and Pablo A. Parrilo\thanks{Amir Ali Ahmadi and Pablo A. Parrilo are with the Laboratory for Information and
Decision Systems, Department of Electrical Engineering and
Computer Science, Massachusetts Institute of Technology. E-mail:
\texttt{a\_a\_a@mit.edu}, \texttt{parrilo@mit.edu}. }  }

\date{}

\begin{document}

\maketitle
\begin{abstract}
A multivariate polynomial $p(x)=p(x_1,\ldots,x_n)$ is sos-convex if
its Hessian $H(x)$ can be factored as $H(x)= M^T(x) M(x)$ with a
possibly nonsquare polynomial matrix $M(x)$. It is easy to see that
sos-convexity is a sufficient condition for convexity of
$p(x)$. Moreover, the problem of deciding sos-convexity of a
polynomial can be cast as the feasibility of a semidefinite program,
which can be solved efficiently.  Motivated by this computational
tractability, it has been recently speculated whether sos-convexity is
also a necessary condition for convexity of polynomials. In this
paper, we give a negative answer to this question by presenting an
explicit example of a trivariate homogeneous polynomial of degree
eight that is convex but not sos-convex. Interestingly, our example is
found with software using sum of squares programming techniques and
the duality theory of semidefinite optimization. As a byproduct of our
numerical procedure, we obtain a simple method for searching over a
restricted family of nonnegative polynomials that are not sums of
squares.
\end{abstract}

\section{Introduction}
In many problems in applied and computational mathematics, we would
like to \emph{decide} whether a multivariate polynomial is convex or
to \emph{parameterize} a family of convex polynomials.  Perhaps the
most obvious instance appears in optimization. It is well known that
in the absence of convexity, global minimization of polynomials is
generally NP-hard \cite{Minimize_poly_Pablo}, \cite{sdprelax},
\cite{nonnegativity_NP_hard}. However, if we somehow know {\it a
priori} that the polynomial is convex, nonexistence of local minima is
guaranteed, and simple gradient descent methods can find a global
minimum. In many other practical settings, we might want to
parameterize a family of convex polynomials that have certain
properties, e.g., that serve as a convex envelope for a non-convex
function, approximate a more complicated function, or fit some data
points with minimum error.  To address many questions of this type, we
need to have an understanding of the algebraic structure of the set of
convex polynomials.

Over a decade ago, Pardalos and Vavasis~\cite{open_complexity} put the
following question proposed by Shor on the list of seven most
important open problems in complexity theory for numerical
optimization: ``Given a degree-$4$ polynomial in $n$ variables, what
is the complexity of determining whether this polynomial describes a
convex function?'' To the best of our knowledge, the question remains
open but the general belief is that the problem should be hard (see
the related work in~\cite{complexity_simplex_convexity}). Not
surprisingly, if testing membership to the set of convex polynomials
is hard, searching and optimizing over them also turns out to be a
hard problem.

The notion of \emph{sos-convexity} has recently been proposed as a
tractable relaxation for convexity based on semidefinite
programming. Broadly speaking, the requirement of positive
semidefiniteness of the Hessian matrix is replaced with the existence
of an appropriately defined sum of squares decomposition. As we will
briefly review in this paper, by drawing some appealing connections
between real algebra and numerical optimization, the latter problem
can be reduced to the feasibility of a semidefinite program. Besides
its computational implications, sos-convexity is an appealing concept
since it bridges the \emph{geometric} and \emph{algebraic} aspects of
convexity. Indeed, while the usual definition of convexity is
concerned only with the geometry of the epigraph, in sos-convexity
this geometric property (or the nonnegativity of the Hessian) must be
\emph{certified} through a ``simple'' algebraic identity, namely the
sum of squares factorization of the Hessian.

Despite the relative recency of the concept of sos-convexity, it has
already appeared in a number of theoretical and practical settings. In
\cite{Helton_Nie_SDP_repres_2}, Helton and Nie use sos-convexity to
give sufficient conditions for semidefinite representability of
semialgebraic sets. In \cite{Lasserre_Jensen_inequality}, Lasserre
uses sos-convexity to extend Jensen's inequality in convex analysis to
linear functionals that are not necessarily probability measures, and
to give sufficient conditions for a polynomial to belong to the
quadratic module generated by a set of polynomials
\cite{Lasserre_Convex_Positive}. More on the practical side, Magnani,
Lall, and Boyd \cite{convex_fitting} have used sum of squares
programming to find sos-convex polynomials that best fit a set of data
points or to find minimum volume convex sets, given by sub-level sets
of sos-convex polynomials, that contain a set of points in space.

Even though it is well-known that sum of squares and nonnegativity are
not equivalent, because of the special structure of the Hessian
matrix, sos-convexity and convexity could potentially turn out to be
equivalent. This speculation has been bolstered by the fact that
finding a counterexample has shown to be difficult and attempts at
giving a non-constructive proof of its existence have seen no success
either. Our contribution in this paper is to give the first such
counterexample, i.e., the first explicit example of a polynomial that
is convex but not sos-convex. This example is presented in
Theorem~\ref{thm:main.result.not_sos_convex}. Our result further
supports the hypothesis that deciding convexity of polynomials should
be a difficult problem. We hope that our counterexample, in a similar
way to what other celebrated counterexamples \cite{MotzkinSOS,
RobinsonSOS, SchmudgenSOS, Choi_Biquadratic} have achieved, will help stimulate
further research and clarify the relationships between the geometric
and algebraic aspects of positivity and convexity.

The organization of the paper is as follows.
Section~\ref{sec:prelim} is devoted to mathematical preliminaries
required for understanding the remainder of this paper. We begin
this section by introducing the cone of nonnegative and sum of
squares polynomials. We briefly discuss the connection between sum
of squares decomposition and semidefinite programming highlighting
also the dual problem. Formal definitions of sos-convex
polynomials and sos-matrices are also given in this section. In
Section~\ref{sec:the.example}, we present our main result, which
is an explicit example of a convex polynomial that is not
sos-convex. Finally, we explain in
Section~\ref{sec:how.we.found.it} how we have utilized
sos-programming together with semidefinite programming duality
theory to find the example presented in
Section~\ref{sec:the.example}. We comment on how one can use the
same methodology to search for a restricted class of nonnegative
polynomials that are not sos.

\section{Background}\label{sec:prelim}
\subsection{Nonnegativity and sum of squares}\label{subsec:nonnegativity.sos}
We denote by $\mathbb{K}[x]\mathrel{\mathop:}=\mathbb{K}[x_1,
\ldots, x_n]$ the ring of polynomials in $n$ variables with
coefficients in the field $\mathbb{K}$. Throughout the paper, we
will have $\mathbb{K}=\mathbb{R}$ or $\mathbb{K}=\mathbb{Q}$. A
polynomial $p(x)\in \mathbb{R}[x]$ is said to be nonnegative or
positive semidefinite (psd) if $p(x)\geq0$ for all
$x\in\mathbb{R}^n$. Clearly, a necessary condition for a
polynomial to be psd is for its total degree to be even. We say
that $p(x)$ is a sum of squares (sos), if there exist polynomials
$q_{1}(x),...,q_{m}(x)$ such that
\begin{equation}\label{eq:sos.decomp.q_i}
p(x)=\sum_{i=1}^{m}q_{i}^{2}(x).
\end{equation}
It is clear that $p(x)$ being sos implies that $p(x)$ is psd. In
1888, David Hilbert \cite{Hilbert_1888} proved that the converse
is true for a polynomial in $n$ variables and of degree $d$
\emph{only} in the following cases:
\begin{itemize}
\item $n=1$ (univariate polynomials of any degree) \item $d=2$
(quadratic polynomials in any number of variables) \item $n=2,\
d=4$ (bivariate quartics)
\end{itemize}
Hilbert showed that in all other cases there exist polynomials
that are psd but not sos. Explicit examples of such polynomials
appeared nearly 80 years later, starting with the celebrated
example of Motzkin, followed by more examples by Robinson, Choi
and Lam, and Lax-Lax and Schm\"{u}dgen. See~\cite{Reznick} for an
outstanding exposition of these counterexamples.

A polynomial $p(x)$ of degree $d$ in $n$ variables has
$l=\binom{n+d}{d}$ coefficients and can therefore be identified with
the $l$-tuple of its coefficients, which we denote by
$\vec{p}\in\mathbb{R}^l$. A polynomial where all the monomials have
the same degree is called a \emph{form}. A form $p(x)$ of degree $d$
is a homogenous function of degree $d$ (since it satisfies $p(\lambda
x)=\lambda^d p(x)$), and has $\binom{n+d-1}{d}$ coefficients. The set
of forms in $n$ variables of degree $d$ is denoted by
$\mathcal{H}_{n,d}$. It is easy to show that if a form of degree $d$
is sos, then $d$ is even, and the polynomials $q_i$ in the sos
decomposition are forms of degree $d/2$. We also denote the set of psd
(resp. sos) forms of degree $d$ in $n$ variables by $P_{n,d}$
(resp. $\Sigma_{n,d}$). Both $P_{n,d}$ and $\Sigma_{n,d}$ are closed
convex cones~\cite{Reznick}, and we have the relation $\Sigma_{n,d}
\subseteq P_{n,d} \subset \mathcal{H}_{n,d}$.

Any form of degree $d$ in $n$ variables can be dehomogenized into
a polynomial of degree $\leq d$ in $n-1$ variables by setting
$x_n=1$. Conversely, any polynomial $p$ of degree $d$ in $n$
variables can be homogenized into a form $p_h$ of degree $d$ in
$n+1$ variables, by adding a new variable $y$, and letting
\[
p_h(x_1,\ldots,x_n,y)\mathrel{\mathop:}=y^{d} \, p\left(\frac{x_1}{y}, \ldots,
\frac{x_n}{y}\right).
\] 
The properties of being psd and sos are preserved under homogenization
and dehomogenization.

An important related problem is Hilbert's 17th problem, which asks if
every psd form must be a sum of squares of rational functions.  In
1927, Artin~\cite{Artin_Hilbert17} answered Hilbert's question in the
affirmative. This result implies that if a polynomial $p(x)$ is psd,
then there must exist an sos polynomial $g(x)$, such that $p(x)g(x)$
is sos. Moreover, Reznick showed in~\cite{Reznick_Unif_denominator}
that if $p(x)$ is \emph{positive definite}, one can always take
$g(x)=(\sum_i x_i^2)^{r}$, for sufficiently large $r$. We will make use
of this key fact in the derivation of our example.

To make the ideas presented so far more concrete, we end this
section by discussing the example of Motzkin. Interestingly
enough, this example will reappear in
Section~\ref{sec:how.we.found.it} where it serves as a starting
point in the numerical procedure that leads to our example.

\begin{example}\label{ex:Motzkin}
The Motzkin polynomial
\begin{equation}\label{eq:Motzkin.poly}
M(x_1,x_2):=x_1^4x_2^2+x_1^2x_2^4-3x_1^2x_2^2+1
\end{equation}
is historically the first known example of a polynomial that is
psd but not sos. Positive semidefiniteness follows from the
arithmetic-geometric inequality, and the nonexistence of an sos
decomposition can be shown by some clever algebraic manipulations
(see~\cite{Reznick} for details). We can homogenize this
polynomial and obtain the Motzkin form
\begin{equation}\label{eq:Motzkin.form}
M_h(x_1,x_2,x_3)\mathrel{\mathop:}=x_1^4x_2^2+x_1^2x_2^4-3x_1^2x_2^2x_3^2+x_3^6,
\end{equation}
which belongs to $P_{3,6} \ \backslash \ \Sigma_{3,6}$ as expected. An
alternative proof of nonnegativity of $M_h$ (resp. $M$) is obtained by
showing that $M_h(x)(x_1^2+x_2^2+x_3^2)$ (resp. $M(x)(x_1^2+x_2^2+1)$)
is sos.  An explicit sos decomposition can be found
in~\cite{PhD:Parrilo}.
\end{example}

In the sequel, we will explain how we can give an alternative
proof of the fact that the Motzkin polynomial is not sos, by
appealing to sos-programming duality.

\subsection{Sum of squares, semidefinite programming, and
duality}\label{subsec:sos.sdp.duality} Deciding nonnegativity of
polynomials is an important problem that arises in many areas of
applied and computational mathematics. Unfortunately, this problem is
known to be NP-hard even when the degree of the polynomial is equal to
four~\cite{sdprelax}, \cite{nonnegativity_NP_hard}. On the other hand,
deciding whether a given polynomial admits an sos decomposition turns
out to be a tractable problem. This tractability stems from the
underlying convexity of the problem as first pointed out in
\cite{Shor}, \cite{NesterovSquared}, \cite{PhD:Parrilo}. More
specifically, it was shown in~\cite{PhD:Parrilo} that one can reduce
the problem of deciding whether a polynomial is sos to feasibility of
a \emph{semidefinite program} (SDP). Semidefinite programs are a
well-studied subclass of convex optimization problems that can be
efficiently solved in polynomial time using interior point
algorithms. Because our space is limited, we refrain from further
discussing SDPs and refer the interested reader to the review
papers~\cite{VaB:96,ToddSDP}. The main theorem that establishes the
link between sum of squares and semidefinite programming is the
following.

\begin{theorem}[\cite{PhD:Parrilo},\cite{sdprelax}]
\label{thm:sos.sdp}
A multivariate polynomial $p(x)$ in $n$ variables and of degree
$2d$ is a sum~of~squares if and only if there exists a positive
semidefinite matrix $Q$ (often called the Gram matrix) such that
\begin{equation}\label{eq:p=z'Qz}
p(x)=z^{T}Qz,
\end{equation}
where $z$ is the vector of monomials of degree up to $d$
\begin{equation*}\label{eq:monomials}
z=[1,x_{1},x_{2},\ldots,x_{n},x_{1}x_{2},\ldots,x_{n}^d].
\end{equation*}
\end{theorem}

Given a polynomial $p(x)$, by expanding the right hand side of
(\ref{eq:p=z'Qz}) and matching coefficients of $p$, we obtain linear
constraints on the entries of $Q$. We also have the constraint that
$Q$ must be a positive semidefinite (PSD\footnote{To avoid potential
confusion, we use the abbreviation psd for positive semidefinite
polynomials and PSD for positive semidefinite matrices. We also denote
a PSD matrix $A$ with the standard notation $A \succeq0$.})
matrix. Therefore, the set described by these constraints is the
intersection of an affine subspace with the cone of PSD matrices. This
is exactly the structure of the feasible set of a semidefinite
program~\cite{VaB:96}. Since the entries of the vector of monomials
$z$ can be algebraically dependent, the matrix $Q$ in the
representation (\ref{eq:p=z'Qz}) is not in general unique. The size of
the matrix $Q$ depends on the size of the vector of monomials. When
there is no sparsity to be exploited $Q$ will have dimensions
$\binom{n+d}{d} \times \binom{n+d}{d}$. If the polynomial $p(x)$ is
homogeneous of degree $2d$, then it suffices to consider in
(\ref{eq:p=z'Qz}) a vector $z$ of monomials of degree exactly
$d$. This will reduce the size of $Q$ to $\binom{n+d-1}{d} \times
\binom{n+d-1}{d}$.

The conversion step of going from an sos decomposition problem to an
SDP problem is fully algorithmic and has been implemented in software
packages such as SOSTOOLS \cite{sostools} and YALMIP
\cite{yalmip}. For instance, we can input a polynomial $p(x)$ into
SOSTOOLS and (if $p(x)$ is sos) it will return a matrix $Q$ and a
vector of monomials $z$. Since $Q$ is PSD, one can compute its
Cholesky factorization $Q=V^{T}V$, which immediately gives rise to an
explicit sos decomposition

\[
p(x)=\sum_{i}(Vz)_{i}^2.
\]

Solutions returned from interior point algorithms are numerical
approximations computed via floating point manipulations. In many
applications in mathematics where the goal is to formally prove a
theorem (as is the case in this paper), it is required to get an exact
algebraic solution. What we mean by this is that given a polynomial
$p(x)\in \mathbb{Q}[x]$, i.e., a polynomial with rational
coefficients, we would like to compute a \emph{rational sos
decomposition}, i.e., a decomposition only consisting of squares of
polynomials in $\mathbb{Q}[x]$. This issue has been studied in detail
in \cite{Rational_SOS_Peyrl_Pablo} where it is shown that the
existence of a rational sos decomposition is equivalent to the
existence of a Gram matrix with rational entries. SOSTOOLS is endowed
with a feature that computes rational decompositions. The work in
\cite{Rational_SOS_Peyrl_Pablo} proposes an efficient mixed
symbolic-numerical approach for this purpose and has been separately
implemented in the package SOS.m2~\cite{SOS.m2} for the computer
algebra system Macaulay~2~\cite{M2}.

Putting the issue of exact computation aside, there are two other
key aspects of sum of squares programming that we would like to
highlight. First, it is not difficult to see that the same
methodology can be used to \emph{search} over sos polynomials in a
convex family of polynomials or even optimize a linear functional
over them~\cite{PhD:Parrilo}. This idea will be crucial in
Section~\ref{sec:how.we.found.it} when we are searching for our
desired polynomial through sos-programming.

The second valuable feature of sos-programming is that when the
semidefinite program arising from Theorem~\ref{thm:sos.sdp} is
infeasible, we get a \emph{certificate} that the polynomial is
\emph{not} sos (though it might still be psd). This certificate is
readily given to us by a feasible solution of the dual semidefinite
program. By definition, the dual cone $\Sigma_{n,d}^{*}$ of the sum of
squares cone $\Sigma_{n,d}$ is the set of all linear functionals $\mu$
that take nonnegative values on it, i.e.,
\begin{equation*}\label{eq:dual_cone_definition}
\Sigma_{n,d}^{*}\mathrel{\mathop:}=\{\mu\in \mathcal{H}_{n,d}^{*}, \ \
\langle \mu,p \rangle \geq 0 \; \; \forall p\in \Sigma_{n,d} \}.
\end{equation*}
Here, the dual space $\mathcal{H}_{n,d}^{*}$ denotes the space of all
linear functionals on $\mathcal{H}_{n,d}$, and $\langle .,. \rangle$
represents the pairing between elements of the primal and the dual
space. If a polynomial is not sos, we can find a dual functional $\mu
\in \Sigma_{n,d}^{*}$ that separates it from the closed convex cone
$\Sigma_{n,d}$. The basic idea behind this is the well known
separating hyperplane theorem in convex analysis; see
e.g.~\cite{BoydBook,Rockafellar}. In Section~\ref{sec:the.example}, we
will see a concrete example of the use of duality when we prove that
our polynomial is not sos-convex. For a more thorough treatment of the
duality theory in semidefinite and sum of squares programming, we
refer to reader to references~\cite{VaB:96} and~\cite{sdprelax},
respectively.

\subsection{Sum of squares matrices and sos-convexity}\label{subsec:sos.convex.sos.matrix}
The notions of positive semidefiniteness and sum of squares of
scalar polynomials can be naturally extended to polynomial
matrices, i.e., matrices with entries in $\mathbb{R}[x]$. We say
that a symmetric polynomial matrix $P(x)\in \mathbb{R}[x]^{m
\times m}$ is PSD if $P(x)$ is PSD for all $x\in \mathbb{R}^n$. It
is straightforward to see that this condition holds if and only if
the polynomial $y^{T}H(x)y$ in $m+n$ variables $[x; y]$ is psd.
The definition of an sos-matrix is as follows
\cite{Kojima_SOS_matrix}, \cite{Symmetry_groups_Gatermann_Pablo},
\cite{matrix_sos_Hol}.

\begin{definition}\label{def:sos-matrix}
A symmetric polynomial matrix $P(x)\in \mathbb{R}[x]^{m \times m},
\ x\in \mathbb{R}^n$ is an sos-matrix if there exists a polynomial
matrix $M(x)\in \mathbb{R}[x]^{s \times m}$ for some
$s\in\mathbb{N}$, such that $P(x)~=~M^{T}(x)M(x)$.
\end{definition}

\begin{lemma}\label{lem:sos-matrix.equivalent.def}
A polynomial matrix $P(x)\in \mathbb{R}[x]^{m \times m}, \ x\in
\mathbb{R}^n$ is an sos-matrix if and only if the scalar
polynomial $y^{T}P(x)y$ is a sum of squares in $\mathbb{R}[x; y]$
.
\end{lemma}

\begin{proof}
One direction is trivial: if $P(x)$ admits the factorization
$M^{T}(x)M(x),$ then the scalar polynomial
$y^{T}M^{T}(x)M(x)y=(M(x)y)^{T}(M(x)y)$ is clearly sos. For the
reverse direction see \cite{Kojima_SOS_matrix}.
\end{proof}
Lemma~\ref{lem:sos-matrix.equivalent.def} enables us to easily check
whether a given polynomial matrix is an sos-matrix with the machinery
explained in Section~\ref{subsec:sos.sdp.duality}. Remarkably, in the
univariate case ($x\in \mathbb{R}$), any PSD polynomial matrix
$P(x)\in \mathbb{R}[x]^{m \times m}$ is an sos-matrix; see
e.g. \cite{CLRrealzeros}. For more details about univariate polynomial
matrices, references to the literature, as well as an efficient
eigenvalue-based method for finding their sos decomposition, we refer
the reader to \cite{SOS_KYP}.

In the multivariate case, however, not every PSD polynomial matrix
must be an sos-matrix. The first counterexample is due to Choi
\cite{Choi_Biquadratic}.  Even though Choi did not have polynomial
matrices in mind, in \cite{Choi_Biquadratic} he showed that not every
psd biquadratic form is a sum of squares of bilinear forms.  His
counterexample can be rewritten as the following polynomial matrix
\begin{equation}\label{eq:Choi.matrix}
C(x)=\begin{bmatrix} x_1^2+2x_2^2&-x_1x_2&-x_1x_3 \\ \\
-x_1x_2&x_2^2+2x_3^2&-x_2x_3 \\ \\
-x_1x_3&-x_2x_3&x_3^2+2x_1^2
\end{bmatrix},
\end{equation}
which is PSD for all $x\in \mathbb{R}^3$ but is not an sos-matrix.

We will now specialize polynomial matrices to Hessians, and discuss
convexity of polynomials. It is well known that a polynomial
$p(x)\mathrel{\mathop:}=p(x_1, \ldots, x_n)$ is convex if and only if
its Hessian
\begin{equation}\label{eq:Hessian}
H(x)=\begin{bmatrix} \frac{\partial^2p}{\partial
x_1^2}&\frac{\partial^2p}{\partial x_1\,\partial
x_2}&\cdots&\frac{\partial^2p}{\partial x_1\,\partial x_n}\\ \\
\frac{\partial^2p}{\partial x_2\,\partial
x_1}&\frac{\partial^2p}{\partial
x_2^2}&\cdots&\frac{\partial^2p}{\partial x_2\,\partial x_n}\\
\\\vdots&\vdots&\ddots&\vdots\\ \\\frac{\partial^2p}{\partial
x_n\,\partial x_1}&\frac{\partial^2p}{\partial x_n\,\partial
x_2}&\cdots&\frac{\partial^2p}{\partial x_n^2}
\end{bmatrix}
\end{equation}
is PSD for all $x\in \mathbb{R}^n$, i.e., is a PSD polynomial
matrix.

\begin{definition}[\cite{Helton_Nie_SDP_repres_2}]
\label{def:sos.convex}
A polynomial $p(x)$ is \emph{sos-convex} if its Hessian $H(x)$ is an sos-matrix.
\end{definition}

Even though we know that not every PSD polynomial matrix is an
sos-matrix, it has been speculated that because of the special
structure of the Hessian as the matrix of the second derivatives,
convexity and sos-convexity of polynomials could perhaps be
equivalent. We will show in the next section that this is not the
case. Note that the example of Choi in (\ref{eq:Choi.matrix}) does
not serve as a counterexample. The polynomial matrix $C(x)$ in
(\ref{eq:Choi.matrix}) is not a valid Hessian, i.e., it cannot be
the matrix of the second derivatives of any polynomial. If this
was the case, the third partial derivatives would commute.
However, we have in particular
$$\frac{\partial C_{1,1}(x)}{\partial x_3}=0\neq-x_3=\frac{\partial C_{1,3}(x)}{\partial x_1}.$$

In \cite{CLRrealzeros}, Choi, Lam, and Reznick generalize the earlier
results of Choi \cite{Choi_Biquadratic} and provide more examples of
psd multiforms that are not sos. Some of their examples can be
rewritten as PSD polynomial matrices that are not sos-matrices. In a
similar fashion, we can show that these matrices also fail to be valid
Hessians.

\section{A polynomial that is convex but not sos-convex}\label{sec:the.example}

We start this section with a lemma that will appear in the proof
of our main result.
\begin{lemma}\label{lem:sos.matrix.then.minor.sos}
If $P(x)\in \mathbb{R}[x]^{m \times m}$ is an sos-matrix, then all
its $2^{m}-1$ principal minors\footnote{We remind the reader that
the principal minors of an $m \times m$ matrix $A$ are the
determinants of all $k \times k$ ($1\leq k \leq m$) sub-blocks
whose rows and columns come from the same index set $S \subset
\{1, \ldots, m\}$.} are sos polynomials. In particular, $\det(P)$
and the diagonal elements of $P$ must be sos polynomials.
\end{lemma}

\begin{proof}
We first prove that $\det(P)$ is sos. By
Definition~\ref{def:sos-matrix}, we have $P(x)=M^{T}(x)M(x)$ for
some $s \times m$ polynomial matrix $M(x)$. If $s=m$, we have
$$\det(P)=\det(M^{T})\det(M)=(\det(M))^2$$ and the result is
immediate. If $s>m$, the result follows from the Cauchy-Binet
formula\footnote{Given matrices $A$ and $B$ of size $m \times s$
and $s \times m$ respectively, the Cauchy-Binet formula states
that
$$\det(AB)=\sum_S \det(A_S)\det(B_S),$$ where $S$ is a subset of $\{1, \ldots, s\}$ with $m$ elements, $A_S$ denotes the $m \times m$ matrix whose columns are the columns of $A$ with index from $S$, and similarly $B_S$ denotes the $m \times m$ matrix whose rows are the rows of $B$ with index from
$S$.}. We have $$ \begin{array}{rll}\det(P)&=&\sum_S \det(M^{T})_S
\det(M_S) \\ \\ \quad&=&\sum_S \det(M_S)^{T} \det(M_S)\\
\\ \quad &=&\sum_S(\det(M_S))^2. \end{array}$$ Finally, when $s<m$, $\det(P)$ is zero
which is trivially sos. In fact, the Cauchy-Binet formula also
holds for $s=m$ and $s<m$, but we have separated these cases for
clarity of presentation.

Next, we need to prove that the smaller minors of $P$ are also
sos. Define $\mathcal{M}=\{1, \ldots, m\},$ and let $I$ and $J$ be
nonempty subsets of $\mathcal{M}$. Denote by $P_{IJ}$ a sub-block
of $P$ with row indices from $I$ and column indices from $J$. It
is easy to see that
$$P_{JJ}=(M^{T})_{J\mathcal{M}}M_{\mathcal{M}J}=(M_{\mathcal{M}J})^TM_{\mathcal{M}J}.$$
Therefore, $P_{JJ}$ is an sos-matrix itself. By the proceeding
argument $\det(P_{JJ})$ must be sos, and hence all the principal
minors are sos.
\end{proof}

\begin{remark}\label{rmk:iff.minor.psd.sos}
The converse of Lemma~\ref{lem:sos.matrix.then.minor.sos} does not
hold. The Choi matrix (\ref{eq:Choi.matrix}) serves as a
counterexample. It is easy to check that all $7$ principal minors of
$C(x)$ are sos polynomials and yet it is not an sos-matrix. This is in
contrast with the fact that a polynomial matrix is PSD if and only if
all its principal minors are psd polynomials. The latter statement
follows almost immediately from the well-known fact that a constant
matrix is PSD if and only if all its principal minors are nonnegative.
\end{remark}

We are now ready to state our main result.
\begin{theorem}\label{thm:main.result.not_sos_convex}
There exists a polynomial that is convex but not sos-convex. In
particular, the trivariate form of degree 8 given by
\begin{equation}\label{eq:our.example.p}
\begin{array}{rlll}
p(x)&=&32x_1^8+118x_1^6x_2^2+40x_1^6x_3^2+25x_1^4x_2^4-43x_1^4x_2^2x_3^2-35x_1^4x_3^4+3x_1^2x_2^4x_3^2
\\
\\ \quad&\ &-16x_1^2x_2^2x_3^4+24x_1^2x_3^6+16x_2^8+44x_2^6x_3^2+70x_2^4x_3^4+60x_2^2x_3^6+30x_3^8
\end{array}
\end{equation}
has these properties.
\end{theorem}
\begin{proof}
Let $H(x)$ denote the Hessian of $p(x)$. Convexity follows from the
fact that
\begin{equation}
(x_1^2+x_2^2+x_3^2) \, H(x)=M^{T}(x)M(x),
\end{equation}
for some polynomial matrix $M(x)$. Equivalently,
\begin{equation}\label{eq:xi^2.y.H.y}
(x_1^2+x_2^2+x_3^2) \, y^{T}H(x)y
\end{equation}
is a sum of squares in $\mathbb{R}[x;y]$, which shows that $H(x)$ is a
PSD polynomial matrix. In the Appendix, we provide an explicit sos
representation in terms of rational Gram matrices for the polynomial
(\ref{eq:xi^2.y.H.y}).  This representation was found using SOSTOOLS
along with the SDP solver SeDuMi~\cite{sedumi}.

To prove that $p(x)$ is not sos-convex, by
Lemma~\ref{lem:sos.matrix.then.minor.sos} it suffices to show that
\begin{equation}
\begin{aligned}
H_{1,1}(x)=\frac{\partial^2 p}{\partial x_1^2} 
&=1792x_1^6+3540x_1^4x_2^2+1200x_1^4x_3^2+300x_1^2x_2^4 \\ 
& \qquad -516x_1^2x_2^2x_3^2-420x_1^2x_3^4+6x_2^4x_3^2-32x_2^2x_3^4+48x_3^6
\end{aligned}
\label{eq:main.ex.H(1,1)}
\end{equation}
is not sos (though it must be psd because of convexity). Define the
subspace $\mathcal{S} \subset \mathcal{H}_{3,6}$ as
\begin{equation}
\mathcal{S}\mathrel{\mathop:}=\mbox{span}\{x_1^6, x_1^4x_2^2,
x_1^4x_3^2, x_1^2x_2^4, x_1^2x_2^2x_3^2, x_1^2x_3^4, x_2^4x_3^2,
x_2^2x_3^4, x_3^6 \},
\label{eq:set.of.monomials.of.H11}
\end{equation}
which are the trivariate sextic forms containing only the
monomials in (\ref{eq:set.of.monomials.of.H11}). Note that $H_{1,1}$
belongs to $\mathcal{S}$. We will prove that $H_{1,1}$ is not sos by
presenting a dual functional $\xi$ that separates $H_{1,1}$ from
$\Sigma_{3,6}\cap\mathcal{S}$.

Consider the vector of coefficients\footnote{As a
trivariate form of degree $6$, $H_{1,1}$ should have $28$
coefficients. We refrain from showing the coefficients that are
zero since our analysis is done in the lower dimensional subspace
$\mathcal{S}$.}  of $H_{1,1}$ with the ordering as written in
(\ref{eq:main.ex.H(1,1)}):
\begin{equation*}\label{eq:vec.of.H(1,1)}
\vec{H}_{1,1}^{T}=[1792, 3540, 1200, 300, -516, -420, 6, -32, 48].
\end{equation*}
Using the same ordering, we can represent our dual functional
$\xi$ with the vector
\begin{equation*}\label{eq:vec.c}
c^{T}=[ 0.039, 0.051, 0.155, 0.839, 0.990, 1.451, 35.488, 20.014,
17.723],
\end{equation*}
which will serve as a separating hyperplane. We have
\begin{equation}\label{eq:c'.H(1,1)=-8.94}
\langle\xi,H_{1,1}\rangle=c^{T}\vec{H}_{1,1}=-8.948<0.
\end{equation}
On the other hand, we claim that for any form
$w\in\Sigma_{3,6}\cap\mathcal{S}$, we should have
\begin{equation}\label{eq:c.w>=0}
\langle\xi,w\rangle=c^{T}\vec{w}\geq0.
\end{equation}
Indeed, if $w$ is sos, by Theorem~\ref{thm:sos.sdp} it can be
written in the form
\[
w(x)=z^{T}Qz= \mathrm{Tr} \ Q \cdot zz^{T},
\]
for some $Q\succeq0$, and a vector of monomials
\[
z^{T}=[x_1^3,  x_1x_2^2,  x_1x_3^2,  x_1^2x_2, x_2x_3^2,
x_1x_2x_3, x_3x_1^2, x_3x_2^2, x_3^3]
\]
that includes all monomials of degree $3$ except for $x_2^3$,
which is not required. It is not difficult to see that
\begin{equation}\label{eq:c.vec(w)=traceQzzz'}
c^{T}\vec{w}= \mathrm{Tr}  \, Q \cdot (zz^{T}) \vert_c,
\end{equation}
where by $(zz^{T})\vert_c$ we mean a matrix where each monomial in
$zz^{T}$ is replaced with the corresponding element of the vector
$c$ (or zero, if the monomial is not in $\mathcal{S}$). This yields the matrix
\[
(zz^{T})\vert_c=
\begin{bmatrix}
    0.039  &  0.051 &   0.155 &        0 &        0  &       0 &        0&         0 &
0\\
    0.051  &  0.839 &   0.990 &        0 &        0  &       0 &        0&         0 &
0\\
    0.155  &  0.990 &   1.451 &        0 &        0  &       0 &        0&         0 &
0\\
         0  &       0 &        0 &   0.051 &        0.990  &  0 &        0&         0 &
0\\
         0  &       0 &        0 &        0.990 &   20.014  &       0 &        0&         0 &
0\\
         0  &       0 &        0 &   0 &        0  & 0.990 &        0&         0 &
0\\
         0  &       0 &        0 &        0 &        0  &       0 &   0.155&    0.990 &
1.451\\
         0  &       0 &        0 &        0 &        0  &       0 &   0.990&   35.488 &
20.014\\
         0  &       0 &        0 &        0 &        0  &       0 &   1.451&   20.014 &  17.723

\end{bmatrix},
\]
and we can easily check that it is positive definite. Therefore,
equation (\ref{eq:c.vec(w)=traceQzzz'}) along with the fact that $Q$
is positive semidefinite implies that (\ref{eq:c.w>=0}) holds. This
completes the proof.
\end{proof}

We end this section with a few remarks on some of the properties
of the polynomial $p(x)$ in~(\ref{eq:our.example.p}).
\begin{remark}
The Gram matrices in the sos decomposition of (\ref{eq:xi^2.y.H.y})
(presented in the Appendix) are positive definite. This shows that for
all nonzero $x$, the Hessian $H(x)$ is positive definite and hence
$p(x)$ is in fact \emph{strictly~convex}; i.e.,
\[
p(\lambda x+(1-\lambda)y)<\lambda p(x)+(1-\lambda) p(y) \qquad
\forall x,y \in \mathbb{R}^3 \quad \mbox{and} \quad \lambda \in (0,1).
\]
\end{remark}

\begin{remark}
Because of strict convexity and the fact that $H_{1,1}$ is
\emph{strictly} separated from $\Sigma_{3,6}$ (see
(\ref{eq:c'.H(1,1)=-8.94})), it follows that $p(x)$ is in the interior
of the set of trivariate forms of degree $8$ that are convex but not
sos-convex. In other words, there exists a neighborhood of polynomials
around $p(x)$, such that \emph{every} polynomial in this neighborhood
is also convex but not sos-convex.

\end{remark}

\begin{remark}
As explained in Section~\ref{subsec:nonnegativity.sos}, we can
dehomogenize the form in (\ref{eq:our.example.p}) into a
polynomial in two variables by letting
\begin{equation}\label{eq:our.example.p.dehomog}
p_{dh}(x_1,x_2)\mathrel{\mathop:}=p(x_1,x_2,1).
\end{equation}
The bivariate polynomial $p_{dh}$ has degree $8$ and we can check
that it is still convex but not sos-convex. It is interesting to
note that $p_{dh}$ is an example with the minimum possible number
of variables since we know that all convex univariate polynomials
are sos-convex. As for minimality in the degree, we do not know if
an example with lower degree exists. However, we should note that
a bivariate form of degree $4$ cannot be convex but not
sos-convex. The reason is that the entries of the Hessian of such
polynomial would be bivariate quadratic forms. It is known that a
matrix with such entries is PSD if and only if it is an
sos-matrix~\cite{Choi_Biquadratic}.

\end{remark}

\begin{remark}
Unlike nonnegativity and sum of squares, sos-convexity may not be
preserved under homogenization. To give a concrete example, one
can check that
$\bar{p}_{dh}(x_2,x_3)\mathrel{\mathop:}=p(1,x_2,x_3)$ is
sos-convex, i.e., the $2 \times 2$ Hessian of
$\bar{p}_{dh}(x_2,x_3)$ is an sos-matrix.
\end{remark}

\begin{remark}
It is easy to argue that the polynomial $p$ in
(\ref{eq:our.example.p}) must itself be nonnegative. Since $p$ is
strictly convex, it has a unique global minimum. Clearly, the
gradient of $p$ has no constant terms and hence vanishes at the
origin. Therefore, $x=0$ must be the unique global minimum of $p$.
Because we have $p(0)=0$, it follows that $p$ is in fact positive
definite.
\end{remark}

\begin{remark}
In \cite{Helton_Nie_SDP_repres_2}, Helton and Nie prove that if a
nonnegative polynomial is sos-convex, then it must be sos. Since
$p$ is not sos-convex, we cannot directly use their result to
claim that $p$ is sos. However, we have independently checked that
this is the case simply by getting an explicit sos decomposition
of $p$ using SOSTOOLS.
\end{remark}

\section{Our procedure for finding the example}
\label{sec:how.we.found.it}

As we mentioned in Section~\ref{subsec:sos.sdp.duality}, one of
the main strengths of sos-programming is in its ability to
\emph{search} over sos polynomials in a convex family of
polynomials. Our main example in~(\ref{eq:our.example.p}) has in
fact been found by solving an sos-program. In this section, we
explain how this has been exactly done.

The task of finding a polynomial $p(x)$ that is convex but not
sos-convex is equivalent to finding a polynomial matrix $H(x)$
that is a \emph{valid Hessian} (i.e., it is a matrix of second
derivatives), and satisfies the following requirement on the
scalar polynomial $y^{T}H(x)y$ in $[x;y]$:
\begin{equation}\label{eq:y'H(x)y}
y^{T}H(x)y \quad \mbox{must be psd but not sos.}
\end{equation}
Indeed, if such a matrix $H(x)$ is found, the desired polynomial
$p(x)$ can be recovered from it by integration. Unfortunately, a
constraint of type (\ref{eq:y'H(x)y}) that requires a polynomial to be
psd but \emph{not} sos cannot be easily handled, since it is a
non-convex constraint. This is easy to see from a geometric viewpoint,
since as Theorem~\ref{thm:sos.sdp} suggests, an sos-program can be
converted to an equivalent semidefinite program. We know that the
feasible set of a semidefinite program is always a convex set. On the
other hand, for a fixed degree and dimension, the set of psd
polynomials that are not sos is generally non-convex. Nevertheless, we
are going to see that by making use of dual functionals of the sos
cone along with Reznick's result on Hilbert's 17th
problem~\cite{Reznick_Unif_denominator}, we can formulate an
sos-program that searches over a \emph{convex subset} of the set of
polynomials that are psd but not sos. The idea behind our algorithm
closely resembles the proof of
Theorem~\ref{thm:main.result.not_sos_convex}.

\paragraph{The algorithm.} The sos-program which has led to our main
result in (\ref{eq:our.example.p}) can be written in pseudo-code as
follows.

\floatname{algorithm}{SOS-Program}
\begin{algorithm*}
\caption{\ } \label{alg:how.we.found.it}
\begin{algorithmic}[1]
\State Parameterize $p(x)$ as a form of degree $8$ in $3$ variables. 
\State Compute the Hessian $H(x)=\frac{\partial^2p}{\partial x^2}$.
\State Impose the constraint
 \begin{equation}\label{eq:algo_xi^2.y.H.y} (x_1^2+x_2^2+x_3^2)^ry^{T}H(x)y \quad \mbox{sos}. \end{equation}
for some integer $r \geq 1$.
\State Impose the constraint
\begin{equation}\label{eq:b.vec(w)<0}\quad \quad \quad \quad \quad \ \ \ \ \ \   \langle \mu, H_{1,1} \rangle <0.\end{equation}
for some (carefully chosen) dual functional $\mu \in \Sigma_{3,6}^*$.
\end{algorithmic}
\end{algorithm*}

The decision variables of this sos-program are the coefficients of
the polynomial $p(x)$ that also appear in the entries of the
Hessian matrix $H(x)$. The scalar $r$ and the dual functional
$\mu$ must be fixed \emph{a priori} as explained in the sequel.
Note that the constraints (\ref{eq:algo_xi^2.y.H.y}) and
(\ref{eq:b.vec(w)<0}) are linear in the decision variables and
indeed the feasible set described by these constraints is a convex
set.

We claim that if this sos-program is feasible, the solution $p(x)$
will be convex but not sos-convex. The requirement of $H(x)$ being
a valid Hessian is met by construction since $H(x)$ is obtained by
twice differentiating a polynomial. It is also easy to see that if
the constraint in (\ref{eq:algo_xi^2.y.H.y}) is satisfied, then
$y^{T}H(x)y$ will be psd. The same implication would hold if instead of
$(x_1^2+x_2^2+x_3^2)^r$ we used any other positive definite
polynomial. As discussed in
Section~\ref{subsec:nonnegativity.sos}, the reason for using this
particular form is due to the result of
Reznick~\cite{Reznick_Unif_denominator}, which states that if
$y^{T}H(x)y$ is positive definite, then (\ref{eq:algo_xi^2.y.H.y})
must be satisfied for sufficiently large $r$. In our case, it was
sufficient to take $r=1$.

In order to guarantee that $y^{T}H(x)y$ is not sos, by
Lemma~\ref{lem:sos.matrix.then.minor.sos} it suffices to require
at least one of the principal minors of $H(x)$ not to be sos
(though they must all be psd because of (\ref{eq:algo_xi^2.y.H.y})
and Remark~\ref{rmk:iff.minor.psd.sos}). The constraint in
(\ref{eq:b.vec(w)<0}) is imposing this requirement on the first
diagonal element $H_{1,1}$. Since $p(x)$ is a form of degree $8$
in $3$ variables, $H_{1,1}$ will be a form of degree $6$ in $3$
variables. The role of the dual functional
$\mu\in\Sigma_{3,6}^{*}$ in (\ref{eq:b.vec(w)<0}) is to separate
$H_{1,1}$ from $\Sigma_{3,6}$. Once an ordering on the monomials
of $H_{1,1}$ is fixed, the inequality in (\ref{eq:b.vec(w)<0}) can
be written as
\begin{equation}
\langle \mu, H_{1,1} \rangle=b^{T}\vec{H}_{1,1}<0,
\end{equation}
where $b\in \mathbb{R}^{28}$ represents our separating hyperplane and
must be fixed \emph{a priori}. We explain next our specific choice of
the dual functional $\mu$.

\paragraph{Finding a separating hyperplane.} There are several ways to
obtain a separating hyperplane for $\Sigma_{3,6}$. In particular,
we can find a dual functional that separates the Motzkin form in
(\ref{eq:Motzkin.form}) from $\Sigma_{3,6}$. This can be done in
at least a couple of different ways. For example, we can formulate
a semidefinite program that requires the Motzkin form to be sos.
This program is clearly infeasible. A feasible solution to its
dual semidefinite program will give us the desired separating
hyperplane. Most SDP solvers, such as SeDuMi, use primal-dual
interior point algorithms to solve an SDP. Therefore, once the
primal SDP is infeasible, a dual feasible solution can readily be
recovered from the solver.

Another way to obtain a separating hyperplane for the Motzkin form
$M_h(x)$ is to find its (Euclidean) projection $M_h^p(x)$ onto the
cone $\Sigma_{3,6}$. Since the projection is done onto a convex set,
the hyperplane tangent to $\Sigma_{3,6}$ at $M_h^p(x)$ will be
supporting $\Sigma_{3,6}$. 
The projection $M_h^p(x)$ can be obtained by searching for an sos
polynomial that is closest in the $2$-norm of the coefficients to
the Motzkin form. This search can be formulated as the following
sos-program:

\begin{equation} \label{eq:sos.prog.projection.of.motzkin}
\begin{array}{lr}
\min & ||\vec{M_h}(x)-\vec{q}(x)||_{2} \\
\mbox{subject to}& q(x) \quad \mbox{sos}.
\end{array}
\end{equation}
Here, $q(x)$ is parameterized as a degree $6$ form in $3$
variables. The objective function in
(\ref{eq:sos.prog.projection.of.motzkin}) can be converted to a
semidefinite constraint using standard tricks; see
e.g.~\cite{VaB:96}.

We have used SOSTOOLS and SeDuMi to obtain a feasible solution to
SOS-Program~\ref{alg:how.we.found.it} with $r=1$ and the dual function
$\mu$ computed using the projection approach described above. In order
to end up with the form in (\ref{eq:our.example.p}) which has integer
coefficients, some post-processing has been done on this feasible
solution. This procedure includes truncation of the coefficients and
some linear coordinate transformations.

We shall end our discussion with a couple of remarks.

\begin{remark}
If the constraints and the objective function of a semidefinite
program possess some type of symmetry, the same symmetry will
generally be inherited in the solution returned by interior point
algorithms. For example, consider the sos-program in
(\ref{eq:sos.prog.projection.of.motzkin}). The Motzkin form
$M_h(x)$ is symmetric in $x_1$ and $x_2$; see
(\ref{eq:Motzkin.form}). Therefore, it turns out that the optimal
solution $M_h^p(x)$ is also symmetric in $x_1$ and $x_2$. On the
other hand, our main example $p(x)$ in (\ref{eq:our.example.p})
and its dehomogenized version $p_{dh}(x)$ in
(\ref{eq:our.example.p.dehomog}) are not symmetric in $x_1$ and
$x_2$. Even though constraint (\ref{eq:algo_xi^2.y.H.y}) and the
dual functional $b$ possess this symmetry, the symmetry is being
broken by imposing constraint (\ref{eq:b.vec(w)<0}) on the second
partial derivative with respect to $x_1$.
\end{remark}

\begin{remark}
Perhaps of independent interest, the methodology explained in this
section can be employed to search or optimize over a restricted
family of psd polynomials that are not sos using sos-programming.
In particular, we can use this technique to simply find more
instances of such polynomials. In order to impose a constraint
that some polynomial $q(x)$ must belong to $P_{n,d} \ \backslash \
\Sigma_{n,d}$, we can use a dual functional
$\eta\in\Sigma^{*}_{n,d}$ to separate $q(x)$ from $\Sigma_{n,d}$,
and then require $q(x)(\sum_{i=1}^{n} x_i^2)^{r}$ to be sos, so
that $q(x)$ stays in $P_{n,d}$.
\end{remark}

\bibliographystyle{abbrv}
\bibliography{AhmadiParrilo}

\newpage

\appendix 

\section{Rational SOS Decomposition of the Hessian}

In this Appendix we present an explicit SOS decomposition of the
Hessian of the form in Theorem~\ref{thm:main.result.not_sos_convex}.
Since the form (\ref{eq:xi^2.y.H.y}) belongs to $\mathcal{H}_{3,10}
\approx \mathbb{R}^{220}$ (is a trivariate form of degree 10), we have
exploited the Newton polytope of the form, as well as its symmetries
in order to reduce the size of our certificate.  For a description of
general techniques for exploiting structure and symmetries in sum of
squares programs, the reader is referred
to~\cite{SOSstructBook,Symmetry_groups_Gatermann_Pablo}.

For our example, note that only even powers of $x_1$, $x_2$, and $x_3$
appear in the form $p(x)$ in (\ref{eq:our.example.p}). Therefore,
$p(x)$ is invariant under changing the signs of $x_1$, $x_2$, and
$x_3$.  This induces a similar type of symmetry in the form
$(x_1^2+x_2^2+x_3^2)\, y^{T}H(x)y$ in (\ref{eq:xi^2.y.H.y}). Namely,
the form is invariant under the group $Z_2^3$, acting via the
transformations $(x_i, y_i) \mapsto (-x_i, -y_i)$ for $i=1,2,3$. This
symmetry enables us to solve a semidefinite program of considerably
smaller size by splitting the monomials in isotypic components
corresponding to the parity of exponents. More concretely, we can
split the monomials of degree five into four groups $z_1,z_2,z_3$, and
$z_4$ that are associated with the exponent patterns (OEE), (EOE),
(EEO), and (OOO). Here, E/O means even/odd, and corresponds to the
parity of the monomials in terms of their \emph{combined} degree in
$x$ and $y$, i.e., the monomials are classified in terms of their
degree in $\{x_1,y_1\}$, $\{x_2,y_2\}$, and $\{x_3,y_3\}$. For
example, the monomial $x_1^2x_2^2y_2$ has combined degrees equal to
$(2,3,0)$ and therefore belongs to the (E,O,E) group.

An explicit sos decomposition of the form in (\ref{eq:xi^2.y.H.y}) is:
\[
(x_1^2+x_2^2+x_3^2) \, y^{T}H(x)y=z_1^{T}Q_1z_1+z_2^{T}Q_2z_2+z_3^{T}Q_3z_3+z_4^{T}Q_4z_4,
\]
where the matrices $Q_i$ and monomial vectors $z_i$ are given by
\[
z_1^{T}= [x_3^4y_1, x_2^2x_3^2y_1, x_1x_3^3y_3, x_1x_2x_3^2y_2,
x_1x_2^2x_3y_3, x_1x_2^3y_2, x_1^2x_3^2y_1, x_1^2x_2^2y_1,
x_1^3x_3y_3, x_1^3x_2y_2, x_1^4y_1],
\]

\setcounter{MaxMatrixCols}{20}
\[
Q_1=
\begin{pmatrix}
                 48&    \frac{-100}{9}&     \frac{161}{3}&     \frac{-53}{3}&    \frac{515}{18}&   \frac{-409}{12}&   \frac{-2134}{9}&  \frac{-8396}{27}&    \frac{-160}{3}&    \frac{-448}{3}&   \frac{-1027}{5}\\
     \frac{-100}{9}&      \frac{32}{9}&   \frac{-217}{18}&     \frac{55}{12}&     \frac{-16}{3}&                 6&   \frac{1549}{27}&     \frac{375}{7}&    \frac{329}{30}&      \frac{85}{3}&   \frac{1079}{27}\\
      \frac{161}{3}&   \frac{-217}{18}&    \frac{3196}{3}&     \frac{260}{3}&     \frac{329}{3}&   \frac{1939}{30}&              -155& \frac{-12281}{30}&  \frac{-3389}{15}&   \frac{-1365}{8}&  \frac{-2069}{12}\\
      \frac{-53}{3}&     \frac{55}{12}&     \frac{260}{3}&    \frac{6862}{9}&    \frac{543}{10}&     \frac{118}{3}&      \frac{32}{3}&     \frac{454}{3}&   \frac{-223}{24}&  \frac{-1117}{27}&     \frac{350}{3}\\
     \frac{515}{18}&     \frac{-16}{3}&     \frac{329}{3}&    \frac{543}{10}&  \frac{16802}{27}&     \frac{434}{3}&  \frac{-3911}{30}&              -182&  \frac{-5185}{27}&    \frac{-341}{3}&              -180\\
    \frac{-409}{12}&                 6&   \frac{1939}{30}&     \frac{118}{3}&     \frac{434}{3}&               768&     \frac{359}{3}&     \frac{866}{3}&                53&               -97&    \frac{1059}{4}\\
    \frac{-2134}{9}&   \frac{1549}{27}&              -155&      \frac{32}{3}&  \frac{-3911}{30}&     \frac{359}{3}&  \frac{22144}{15}&  \frac{31796}{27}&     \frac{965}{4}&    \frac{2042}{3}&    \frac{2815}{3}\\
   \frac{-8396}{27}&     \frac{375}{7}& \frac{-12281}{30}&     \frac{454}{3}&              -182&     \frac{866}{3}&  \frac{31796}{27}&   \frac{24706}{9}&    \frac{1085}{3}&  \frac{13541}{12}&    \frac{4441}{3}\\
     \frac{-160}{3}&    \frac{329}{30}&  \frac{-3389}{15}&   \frac{-223}{24}&  \frac{-5185}{27}&                53&     \frac{965}{4}&    \frac{1085}{3}&    \frac{1682}{9}&   \frac{1861}{12}&    \frac{1156}{3}\\
     \frac{-448}{3}&      \frac{85}{3}&   \frac{-1365}{8}&  \frac{-1117}{27}&    \frac{-341}{3}&               -97&    \frac{2042}{3}&  \frac{13541}{12}&   \frac{1861}{12}&    \frac{8158}{9}&    \frac{3044}{3}\\
    \frac{-1027}{5}&   \frac{1079}{27}&  \frac{-2069}{12}&     \frac{350}{3}&              -180&    \frac{1059}{4}&    \frac{2815}{3}&    \frac{4441}{3}&    \frac{1156}{3}&    \frac{3044}{3}&      1792

\end{pmatrix},
\]

\[
z_2^{T}= [x_3^4y_2,
    x_2x_3^3y_3,
    x_2^2x_3^2y_2,
    x_2^3x_3y_3,
    x_2^4y_2,
    x_1x_2x_3^2y_1,
    x_1x_2^3y_1,
    x_1^2x_3^2y_2,
    x_1^2x_2x_3y_3,
    x_1^2x_2^2y_2,
    x_1^3x_2y_1,
    x_1^4y_2],
\]

\[
Q_2=
\begin{pmatrix}
        120&        86&        21&        36&     \frac{312}{5}&    \frac{-206}{3}&    \frac{763}{12}&    \frac{-341}{3}&       \frac{5}{3}&     \frac{550}{9}&     \frac{200}{3}&      \frac{14}{3}\\
         86&   \frac{13252}{9}&    \frac{1459}{3}&       334&   \frac{5293}{12}&   \frac{2369}{18}&       -76&        45&       -76&   \frac{5489}{30}&  \frac{-5501}{30}&  \frac{-4991}{24}\\
         21&    \frac{1459}{3}&  \frac{18586}{15}&    \frac{1427}{4}&    \frac{3442}{9}&  \frac{-1777}{12}&     \frac{320}{3}&     \frac{349}{9}&   \frac{-851}{30}&     \frac{479}{3}&     \frac{332}{3}&  \frac{-3169}{27}\\
         36&       334&    \frac{1427}{4}&    \frac{7220}{9}&       427&        -7&      \frac{70}{3}&   \frac{1189}{30}&  \frac{-4672}{27}&       185&      \frac{91}{3}&    \frac{-257}{3}\\
      \frac{312}{5}&   \frac{5293}{12}&    \frac{3442}{9}&       427&       896&     \frac{-62}{3}&         0&     \frac{239}{3}&      \frac{91}{3}&        64&     \frac{-35}{3}&    \frac{-460}{3}\\
     \frac{-206}{3}&   \frac{2369}{18}&  \frac{-1777}{12}&        -7&     \frac{-62}{3}&  \frac{10000}{27}&  \frac{-6841}{21}&        91&    \frac{329}{30}&      -134& \frac{-18307}{27}&    \frac{-328}{3}\\
     \frac{763}{12}&       -76&     \frac{320}{3}&      \frac{70}{3}&         0&  \frac{-6841}{21}&       300&    \frac{-245}{3}&       -19&       123&    \frac{4927}{9}&     \frac{351}{4}\\
     \frac{-341}{3}&        45&     \frac{349}{9}&   \frac{1189}{30}&     \frac{239}{3}&        91&    \frac{-245}{3}&    \frac{1300}{3}&   \frac{1217}{24}&   \frac{1493}{27}&    \frac{-571}{3}&    \frac{-440}{3}\\
        \frac{5}{3}&       -76&   \frac{-851}{30}&  \frac{-4672}{27}&      \frac{91}{3}&    \frac{329}{30}&       -19&   \frac{1217}{24}&   \frac{8306}{27}&        70&    \frac{-253}{3}&    \frac{-573}{4}\\
      \frac{550}{9}&   \frac{5489}{30}&     \frac{479}{3}&       185&        64&      -134&       123&   \frac{1493}{27}&        70&    \frac{2402}{3}&   \frac{4021}{12}&   \frac{-1667}{9}\\
      \frac{200}{3}&  \frac{-5501}{30}&     \frac{332}{3}&      \frac{91}{3}&     \frac{-35}{3}& \frac{-18307}{27}&    \frac{4927}{9}&    \frac{-571}{3}&    \frac{-253}{3}&   \frac{4021}{12}&    \frac{7114}{3}&    \frac{1204}{3}\\
       \frac{14}{3}&  \frac{-4991}{24}&  \frac{-3169}{27}&    \frac{-257}{3}&    \frac{-460}{3}&    \frac{-328}{3}&     \frac{351}{4}&    \frac{-440}{3}&    \frac{-573}{4}&   \frac{-1667}{9}&    \frac{1204}{3}&       236
\end{pmatrix},
\]

\[
z_3=[x_3^4y_3,
    x_2x_3^3y_2,
    x_2^2x_3^2y_3,
    x_2^3x_3y_2,
    x_2^4y_3,
    x_1x_3^3y_1,
    x_1x_2^2x_3y_1,
    x_1^2x_3^2y_3,
    x_1^2x_2x_3y_2,
    x_1^2x_2^2y_3,
    x_1^3x_3y_1,
    x_1^4y_3],
\]

\[
Q_3=
\begin{pmatrix}
      1680&      634&   \frac{9034}{9}&   \frac{2245}{3}&    \frac{317}{3}&    \frac{703}{3}& \frac{-2671}{18}&   \frac{2002}{3}&    \frac{295}{3}&    \frac{209}{3}&     -166&  \frac{-1113}{5}\\
       634&      918&   \frac{1708}{3}&  \frac{5971}{15}&   \frac{613}{12}&        0&   \frac{155}{12}&    \frac{968}{3}&   \frac{-593}{9}&   \frac{313}{10}&     -149&   \frac{-893}{8}\\
    \frac{9034}{9}&   \frac{1708}{3}&   \frac{5282}{3}&  \frac{9589}{12}&    \frac{566}{9}&  \frac{2537}{18}&   \frac{-146}{3}&   \frac{1787}{3}&  \frac{3619}{30}&   \frac{611}{27}& \frac{-4411}{30}& \frac{-5608}{27}\\
    \frac{2245}{3}&  \frac{5971}{15}&  \frac{9589}{12}&  \frac{13060}{9}&      101&   \frac{863}{12}&   \frac{-208}{3}& \frac{14299}{30}&    \frac{-68}{3}&       59&     -124&   \frac{-490}{3}\\
     \frac{317}{3}&   \frac{613}{12}&    \frac{566}{9}&      101&       88&     \frac{37}{3}&    \frac{-16}{3}&  \frac{2402}{27}&       27&   \frac{-181}{3}&    \frac{-59}{3}&    \frac{-41}{3}\\
     \frac{703}{3}&        0&  \frac{2537}{18}&   \frac{863}{12}&     \frac{37}{3}&    \frac{920}{9}& \frac{-2231}{27}&    \frac{307}{3}&     \frac{-7}{3}&   \frac{163}{10}& \frac{-2141}{15}&   \frac{-151}{4}\\
  \frac{-2671}{18}&   \frac{155}{12}&   \frac{-146}{3}&   \frac{-208}{3}&    \frac{-16}{3}& \frac{-2231}{27}&  \frac{1762}{21}&  \frac{-627}{10}&     \frac{17}{3}&    \frac{-41}{3}&  \frac{1727}{27}&     \frac{65}{3}\\
    \frac{2002}{3}&    \frac{968}{3}&   \frac{1787}{3}& \frac{14299}{30}&  \frac{2402}{27}&    \frac{307}{3}&  \frac{-627}{10}& \frac{17956}{15}&  \frac{2353}{24}&  \frac{1988}{27}& \frac{-1333}{12}&  \frac{-2371}{9}\\
     \frac{295}{3}&   \frac{-593}{9}&  \frac{3619}{30}&    \frac{-68}{3}&       27&     \frac{-7}{3}&     \frac{17}{3}&  \frac{2353}{24}& \frac{18134}{27}&     \frac{91}{3}&      182&   \frac{-237}{4}\\
     \frac{209}{3}&   \frac{313}{10}&   \frac{611}{27}&       59&   \frac{-181}{3}&   \frac{163}{10}&    \frac{-41}{3}&  \frac{1988}{27}&     \frac{91}{3}&    \frac{586}{3}&   \frac{-229}{3}&   \frac{-178}{3}\\
      -166&     -149& \frac{-4411}{30}&     -124&    \frac{-59}{3}& \frac{-2141}{15}&  \frac{1727}{27}& \frac{-1333}{12}&      182&   \frac{-229}{3}&   \frac{3346}{3}&    \frac{284}{3}\\
   \frac{-1113}{5}&   \frac{-893}{8}& \frac{-5608}{27}&   \frac{-490}{3}&    \frac{-41}{3}&   \frac{-151}{4}&     \frac{65}{3}&  \frac{-2371}{9}&   \frac{-237}{4}&   \frac{-178}{3}&    \frac{284}{3}&       80
\end{pmatrix},
\]

\[
z_4= [x_2x_3^3y_1,
    x_2^3x_3y_1,
    x_1x_3^3y_2,
    x_1x_2x_3^2y_3,
    x_1x_2^2x_3y_2,
    x_1x_2^3y_3,
    x_1^2x_2x_3y_1,
    x_1^3x_3y_2,
    x_1^3x_2y_3],
\]
and
\[
Q_4=
\begin{pmatrix}
     \frac{344}{9}&   \frac{-133}{9}&     \frac{67}{3}&   \frac{347}{18}&    \frac{-11}{4}&       16& \frac{-8990}{27}&   \frac{-118}{3}&  \frac{-371}{30}\\
    \frac{-133}{9}&        6&   \frac{-97}{12}&    \frac{-22}{3}&      \frac{4}{3}&       -6&  \frac{2630}{21}&     \frac{44}{3}&     \frac{14}{3}\\
      \frac{67}{3}&   \frac{-97}{12}&    \frac{946}{3}&    \frac{113}{3}&     \frac{61}{9}& \frac{-2191}{30}&   \frac{-709}{3}&   \frac{-841}{3}&  \frac{1673}{24}\\
    \frac{347}{18}&    \frac{-22}{3}&    \frac{113}{3}&      930&  \frac{1473}{10}&  \frac{2087}{27}& \frac{-3751}{30}&  \frac{-455}{24}& \frac{-4771}{27}\\
     \frac{-11}{4}&      \frac{4}{3}&     \frac{61}{9}&  \frac{1473}{10}&      844&      106&    \frac{-32}{3}& \frac{-2899}{27}&   \frac{-200}{3}\\
        16&       -6& \frac{-2191}{30}&  \frac{2087}{27}&      106&    \frac{644}{3}&   \frac{-398}{3}&      128&     -124\\
  \frac{-8990}{27}&  \frac{2630}{21}&   \frac{-709}{3}& \frac{-3751}{30}&    \frac{-32}{3}&   \frac{-398}{3}& \frac{81458}{27}&   \frac{1177}{3}&    \frac{280}{3}\\
    \frac{-118}{3}&     \frac{44}{3}&   \frac{-841}{3}&  \frac{-455}{24}& \frac{-2899}{27}&      128&   \frac{1177}{3}&   \frac{1330}{3}& \frac{-1495}{12}\\
   \frac{-371}{30}&     \frac{14}{3}&  \frac{1673}{24}& \frac{-4771}{27}&   \frac{-200}{3}&     -124&    \frac{280}{3}& \frac{-1495}{12}&    \frac{338}{3}
\end{pmatrix}.
\]

\end{document}